\documentclass[10pt]{amsart}
\usepackage{eufrak}
\usepackage{amsfonts}
\usepackage{float}
\usepackage{amsmath,amssymb}
\usepackage{lscape}
\usepackage[dvipdfmx]{graphicx, color}
\usepackage[all]{xy}

\usepackage{comment}

\def\qed{\hfill $\Box$}
\def\proof{\noindent {\sl Proof} :\;  }

\newcommand{\U}{\mathcal{U}}

\newcommand{\A}{\mathcal{A}}

\newcommand{\Lcal}{\mathcal{L}}

\newcommand{\R}{\mathbb{R}}

\newcommand\rank{\mbox{\rm rank}\,}

\def\qed{\hfill $\Box$}
\def\proof{\noindent {\sl Proof} :\;  }

\def\rd{\partial}

\def\bx{\mbox{\boldmath $x$}}

\def\bb{\mbox{\boldmath $b$}}
\def\bc{\mbox{\boldmath $c$}}

\def\bp{\mbox{\boldmath $p$}}

\def\b0{\mbox{\boldmath $0$}}

\def\bbx{\mbox{\tiny$\bx$}}

\def\bbp{\mbox{\tiny$\bp$}}

\newtheorem{thm}{\bf Theorem}[section]
\newtheorem{cor}[thm]{\bf Corollary}

\newtheorem{prop}[thm]{\bf Proposition}
\newtheorem{dfn}[thm]{\bf Definition}

\newtheorem{rem}[thm]{\bf Remark}
\newtheorem{exam}[thm]{\bf Example}

\newcommand{\wt}[1]{\widetilde{#1}}
\newcommand{\vp}{\varphi}

\begin{document}
\title[Local normal forms of em-wavefronts in affine coordinates]
{Local normal forms of em-wavefronts in affine flat coordinates}
\author[N.~Nakajima]{Naomichi Nakajima}
\address[N.~Nakajima]{D2, Graduate School of Information Science and Technology, Hokkaido University,
Sapporo 060-0814, Japan}
\email{nakajima-n@ist.hokudai.ac.jp}
\subjclass[2020]{Primary~57R45; Secondary~53A15, 53B12}
\keywords{wavefronts, singularity theory, affine differential geometry, information geometry, statistical manifold}
\dedicatory{}
\begin{abstract}
In our previous work, we have generalized the notion of dually flat or Hessian manifold to {\em quasi-Hessian manifold}; it admits the Hessian metric to be degenerate but possesses a particular symmetric cubic tensor (generalized Amari-Centsov tensor).  
Indeed, it naturally appears as a singular model in information geometry and related fields. 
A quasi-Hessian manifold is locally accompanied with a possibly multi-valued potential and its dual, whose graphs are called the {\em $e$-wavefront} and the {\em $m$-wavefront} respectively, together with coherent tangent bundles endowed with flat connections.
In the present paper, using those connections and the metric, we give coordinate-free criteria for detecting local diffeomorphic types of $e/m$-wavefronts, and then derive the local normal forms of those (dual) potential functions for the $e/m$-wavefronts in affine flat coordinates by means of Malgrange's division theorem.  
This is motivated by an early work of Ekeland on non-convex optimization and Saji-Umehara-Yamada's work on Riemannian geometry of wavefronts. 
Finally, we reveal a relation of our geometric criteria with information geometric quantities of statistical manifolds. 
\end{abstract}

\maketitle

\section{Introduction}\label{sec1:intro}
In information geometry \cite{AmariNagaoka00}, a dually flat manifold is a smooth manifold endowed with a pseudo-Riemannian metric and two flat affine connections satisfying a certain duality \cite{AmariNagaoka00, Amari16}.
This is also known as a Hessian manifold in affine differential geometry (Shima \cite{Shima}).

As a generalization \cite{NT}, we have introduced a quasi-Hessian manifold to be a smooth manifold $M$ endowed with a possibly degenerate symmetric $(0, 2)$-tensor $h$, called a quasi-Hessian metric, and two coherent tangent bundles $E$ and $E'$ with flat connections $\nabla^E$ and $\nabla^{E'}$, respectively.
In the case where $h$ is non-degenerate everywhere, it coincides with a dually flat (Hessian) manifold, i.e., $TM=E=E'$ and $\nabla^E$ and $\nabla^{E'}$ are mutually dual flat connections on $TM$.
The quasi-Hessian manifold $M$ is locally a Legendre submanifold in the standard contact manifold $\R^{2n+1}=T^*\R^n\times \R$ with the double fibration structure, called a local model. 
The projections of a local model are called $e/m$-Legendre maps, and images of them are called $e/m$-wavefronts, respectively, and coherent tangent bundles $E$ and $E'$ are formed by limiting tangent spaces of the wavefronts.

In the present paper, we characterize typical singularities of these $e/m$-wavefronts in terms of geometric quantities.
Specifically, we give in \S\ref{sec:normal forms} coordinate-free criteria
using $\nabla^E, \nabla^{E'}$ and the metric, for detecting the geometric feature of singularities of $e/m$-wavefronts, and then derive the local normal forms of them in affine flat coordinates.

We mention two preceding related works. 
First, singularities of wavefronts have been investigated by Saji-Umehara-Yamada from the viewpoint of Riemannian geometry.
They gave useful criteria for singularities of types cuspidal edge and swallowtail, and investigated metric geometry of these singularities by introducing several notions of curvatures \cite{KRSUY, SUY}.
We are seeking for an affine differential geometry counterpart and the present paper is the first step in this project.
Secondly, an early work of Ekeland gives a local normal form of the wavefront at an inflectional point in the context of non-convex optimization and calculus of variations \cite{Ekeland}.
In fact, our first normal form (Theorem \ref{m-wavefront1}) refines the result of Ekeland, added some geometric characterization.
Our second normal form (Theorem \ref{m-wavefront2}) is for the case where $h$ is positive semi-definite, e.g., $h$ is a degenerate Fisher-Rao metric in information geometry. 

Statistical manifolds are pseudo-Riemannian manifolds endowed with symmetric cubic tensors, among which dually flat manifolds are special ones   \cite{AmariNagaoka00, Amari16}. 
Namely, the geometry of statistical manifolds is governed by the second and third-order symmetric tensors. 
Fourth-order tensors on a statistical manifold would be also of particular interest, as Eguchi \cite{Eguchi} has explored them thoroughly. In our theorems, the geometric (coordinate-free) criteria for detecting typical singularities of $e/m$-wavefronts are actually concerning the third and fourth derivatives of the canonical divergence,  that is surely related to those studies on third and fourth-order tensors. We describe it explicitly in the final section \S \ref{sec:expression of criteria}.

Throughout the present paper, we assume that manifolds and maps are of class $C^\infty$ for the simplicity and bold letters denote column vectors, e.g., $\bx=(x_1,\cdots,x_n)^T$.

This work was supported by the Hokkaido University Ambitious Doctoral Fellowship (Information Science and AI).

The author thanks his advisor, Professor Toru Ohmoto, for valuable comments and many instructions.

\section{Contact geometry and quasi-Hessian manifolds}\label{sec:2}
In this section, we briefly summarize quasi-Hessian manifold theory in order to describe main results in the next section. 
See \cite{NT} for the detail.

Let $N$ be a $(2n+1)$-dimensional manifold and $\xi$ be a hyperplane field on it; $\xi_p \subset T_pN\;(p\in N)$. $(N, \xi)$ is called a {\em contact manifold} if $\xi$ is locally given by the kernel of $1$-form $\theta$ satisfying $\theta \wedge (d\theta)^n \neq 0$. 
Then, $\xi$ is called a contact structure and $\theta$ is called a (local) contact form. 
For a $(2n+1)$-dimensional contact manifold $(N, \xi)$, a submanifold $L\subset N$ with $\dim L=n$ is called a {\em Legendre submanifold} if $T_pL\subset \xi_p$ for any $p\in L$. 

Let $(\bx,\bp,z)$ denote the standard coordinates of $\R^{2n+1}(=T^*\R^n \times\R)$, where $\bx$ and $\bp$ denote coordinates of the base and fiber space of $T^*\R^n$ respectively. 
We write $\R^n_{\bbx}$ and $\R^n_{\bbp}$ as these spaces respectively in order to distinguish them.
The $1$-form
$$\theta:=dz-\bp^Td\bx = dz-\sum_{i=1}^np_idx_i$$
gives a contact form on $\R^{2n+1}$. Then, $(\R^{2n+1}, \theta)$ is called the standard contact manifold. 

The following fibration is a Legendre fibration (each fiber is a Legendre submanifold):
$$\pi:\R^{2n+1}\to\R^n_{\bbx}\times\R_z,\;(\bx,\bp,z)\mapsto(\bx,z).$$
Consider transformation $\Lcal: \R^{2n+1}\to\R^{2n+1}$ defined by
$$\Lcal(\bx,\bp,z)=(\bx',\bp',z')=(\bp,\bx,\bp^T\bx-z).$$
This $\Lcal$ is a diffeomorphism preserving contact hyper planes; $\Lcal^*\theta=-\theta$. 
The projection along fiber 
$$\pi' :=\pi\circ\Lcal:\R^{2n+1}\to\R^n_{\bbp}\times\R_{z'},\;(\bx,\bp,z)\mapsto(\bp,\bp^T\bx-z)$$
is also a Legendre fibration. Then, the following diagram is called the {\em double fibration structure} for the standard contact manifold:
\begin{align}\label{double fibration}
\xymatrix{\R^{n}_{\bbx}\times \R & \R^{2n+1} \ar[l]_-\pi \ar[r]^-{\pi'} & \R^{n}_{\bbp}\times \R}
\end{align}

Let $L\subset \R^{2n+1}$ be a Legendre submanifold. Legendre maps given by the double fibration structure
$$\pi^e:=\pi\circ \iota: L \to \R^n_{\bbx}\times \R_z, \quad \pi^m:=\pi'\circ \iota:  L \to \R^n_{\bbp}\times \R_{z'}$$
are called the {\em $e/m$-Legendre maps} respectively, where $\iota:L\to\R^{2n+1}$ is the inclusion map.

\begin{dfn}[\cite{NT}]\upshape
For $e/m$-Legendre maps $\pi^e$ and $\pi^m$, 
$$W_e(L):=\pi^e(L)\subset \R^n_{\bbx}\times \R_z, \;\;\; 
W_m(L):=\pi^m(L)\subset \R^n_{\bbp}\times \R_{z'}$$
are called the {\em $e/m$-wavefronts associated to $L$} respectively.
\end{dfn}

The vector bundle $E\,(=E_L)$ on $L$ is defined by
$$E := \{\;(p, w) \in L\times (\R^n_{\bbx} \times \R_z)\; \mid \;dz_p(w)-\bp(p)^T d\bx_p(w)=0\; \}.$$
Since $L$ is a Legendre submanifold, $d\pi^e(T_pL)\subset E_p$ for any $p\in L$, where $E_p$ is the fiber of $p$. 
Therefore, the vector bundle homomorphism
$$\Phi : TL \to E, \quad v_p \mapsto d\pi^e_p(v_p)$$
makes sense.

Let $\wt{\nabla}$ be  
the natural connection on (the trivial vector bundle) $\R^n_{\bbx}\times\R$ over $L$
and $\psi_{p} : \R^n_{\bbx} \times\R \to E_{p}$ be the projection along $z-$axis for $p\in L$.
The connection of $E$, $\nabla^E$, is defined by
$$\nabla^E_X\eta (p) := \psi_{p} \circ  \widetilde{\nabla}_X\eta(p),$$
where $X$ is a vector field on $L$ and $\eta$ is a section of $E$. 
This connection enjoys the following property. 

\begin{prop}[\cite{NT}]\upshape
The connection $\nabla^E$ is flat and `relatively torsion-free', i.e.,  
for any vector fields $X, Y$ on $L$, it holds that 
$$\nabla^E_X(\Phi(Y))-\nabla^E_Y(\Phi(X))=\Phi([X, Y]).$$
\end{prop}

According to Saji-Umehara-Yamada \cite{SUY},  abusing words, we call 
$(E, \Phi, \nabla^E)$ the {\em coherent tangent bundle} associated to $e$-wavefront $W_e(L)$.
In the same way, the coherent tangent bundle $(E', \Phi', \nabla^{E'})$ associated to the $m$-wavefront is defined. 
Note that $\nabla^{E}$ and $\nabla^{E'}$ are flat.

Those vector bundles $E, E'$ are actually defined on $\R^{2n+1}$. Note that the contact hyperplane $\xi$ has a direct decomposition 
$$\xi_p= \ker d\pi'_p\oplus\ker d\pi_p \simeq E_p \oplus E'_p \simeq \R^n_{\bbx}\oplus \R^n_{\bbp}.$$
The  
pseudo-Riemannian metric $\tau$ of type $(n, n)$ on $\xi$ is naturally induced from
\begin{align}\label{tau} 
\tau:=\sum_{i=1}^n dx_idp_i=\frac{1}{2}\sum_{i=1}^n (dx_i \otimes dp_i+dp_i \otimes dx_i).
\end{align}
In fact, there are canonical frames of flat sections for both $E$ and $E'$ given by
$$s_i(p) = \frac{\rd}{\rd x_i} + p_i\frac{\rd}{\rd z} \in E_p, \quad s_i^*(p)= \frac{\rd}{\rd p_i}+x_i\frac{\rd}{\rd z'} \in E'_p,$$
which are projected to $\frac{\rd}{\rd x_i}$, $\frac{\rd}{\rd p_i}$, respectively,  
and satisfy e.g., $\tau(s_i, s_j^*)=\frac{1}{2}\delta_{ij}$.

\begin{dfn}[\cite{NT}]\label{quasi-Hessian metric}\upshape
The {\em quasi-Hessian metric $h$} on $L$ is defined by the pullback of $\tau$:
$$h(Y, Z):=\tau(\iota_*Y, \iota_*Z)\;\;\;  (Y, Z \in TL),$$
where $\iota_*=\Phi\oplus \Phi': TL \hookrightarrow \xi=E\oplus E'$ is the inclusion map.
\end{dfn}

In general, a Legendre submanifold $L\subset \R^{2n+1}$ can be locally expressed by a generating function $g(\bx_I,\bp_J)$\cite{AGV}:
$$\textstyle L=\{(\bx_I,\bx_J,\bp_I,\bp_J,z)\in\R^{2n+1} \mid \bp_I= \frac{\rd g}{\rd \bbx_I},  \bx_J=- \frac{\rd g}{\rd \bbp_J}, z=\bp_J^T\bx_J+g(\bx_I, \bp_J)\},$$
where $I\sqcup J=\{1,\cdots,n\}$ is a partition, $(\bx_I,\bp_J)$ are local coordinates of $L$ and $\frac{\rd g}{\rd \bbx_I}$ denotes the column vector $(\frac{\rd g}{\rd x_i})_{i\in I}^T$.

\begin{prop}[\cite{NT}]\label{quasi-Hessian_local}\upshape
Let $g(\bx_I,\bp_J)$ be a generating function. Then, 
$$h=\sum_{i, k\in I}\frac{\rd^2 g}{\rd x_i\rd x_{k}}\,dx_idx_{k} - \sum_{j, l\in J}\frac{\rd^2 g}{\rd p_j\rd p_{l}}\,dp_jdp_{l}.$$
\end{prop}
Affine transformations on $\R^n_{\bbx}\times\R_z$ and $\R^n_{\bbp}\times\R_{z'}$, respectively, 
$$F(\bx, z)=(A\bx+ \bb,  z+\bc^T\bx+ d),\quad F^*(\bp, z')= (A'\bp+ \bb', z'+\bc'^{T}\bp+ d')$$
determine an {\em affine Legendre equivalence} $\Lcal_F : \R^{2n+1} \to \R^{2n+1}$ defined by
$$\Lcal_F (\bx, \bp, z) = (A\bx+ \bb, \;  A'\bp+ \bb', \; z+\bc^T\bx+ d),$$
where $A$ is an invertible matrix, $A'=(A^T)^{-1}, \; \bb'=A'\bc, \; \bb=A\bc'$, and $d'=\bb'^T\bb-d$.
This transformation $\Lcal_F$ preserves the double fibration structure and the contact structure (and thus $\omega$ and $\tau$).

Suppose two Legendre submanifolds $L_1, L_2$ are identified by $\Lcal_F$; $\Lcal_F(L_1)=L_2$. 
Then, $\Lcal_F$ preserves the quasi-Hessian metric, and induces isomorphisms $E_{L_1} \simeq E_{L_2}$, $E'_{L_1} \simeq E'_{L_2}$ which also identify flat affine connections.
Hence, the following construction makes sense.

\begin{dfn}[\cite{NT}]\upshape
A {\em quasi-Hessian manifold} $(M, \U=\{L_\alpha\})$ is defined by gluing any Legendre submanifolds $L_\alpha \subset \R^{2n+1}$ via affine Legendre equivalences. Then, a possibly degenerate $(0, 2)-$tensor $h$ and coherent tangent bundles $(E, \Phi, \nabla^E)$, $(E', \Phi', \nabla^{E'})$ on $M$ induced from each structure of $L_\alpha$ are well-defined. We call each $L_\alpha$ a local model of $M$.
\end{dfn}
The construction above also defines the contact manifold including a quasi-Hessian manifold $M$. 

Let $(M,h,(E,\Phi, \nabla^E), (E', \Phi', \nabla^{E'}))$ be a quasi-Hessian manifold. For vector fields $X, Y, Z$ on $M$, set
$$(\eta,\eta'):=(\Phi\oplus\Phi')(Y),\;(\zeta,\zeta'):=(\Phi\oplus\Phi')(Z),\;\tau(\eta,\zeta'):=\tau(\eta\oplus 0,0\oplus\zeta').$$

\begin{dfn}[\cite{NT}]\upshape\label{cubic tensor}
For a quasi-Hessian manifold $M$, the {\em canonical cubic tensor} $C$ is defined by
$$C(X, Y, Z):=\tau(\eta, \nabla^{E'}_X\zeta')+\tau(\zeta, \nabla^{E'}_X\eta')
-\tau(\nabla^E_X\eta, \zeta')-\tau(\nabla^E_X\zeta, \eta').$$
\end{dfn}

\begin{prop}[\cite{NT}]\upshape
The canonical cubic tensor $C$ is locally the third derivative of a generating function $g(\bx_I, \bp_J)$: for any $k, l, m$, 
$$C(\rd_k, \rd_l, \rd_m)=\rd_k\rd_l\rd_m g,$$ 
where $\rd_k:=\frac{\rd}{\rd x_k}\; (k\in I) \;\mbox{or}\; \rd_k:=\frac{\rd}{\rd p_k}\; (k \in J)$. Thus $C$ is symmetric.
\end{prop}

\section{Local normal forms of $e/m$-wavefronts}\label{sec:normal forms}
In this section, we characterize the singularities of $e/m$-wavefronts in terms of the quasi-Hessian geometric quantities.  
In particular, we present local normal forms of  (dual) potential functions for $e/m$-wavefronts in affine coordinates. 
Note that if the wavefront is singular, the potential and the dual potential may be multi-valued. 
A major tool is  Malgrange's division theorem in singularity theory.

\begin{thm}[Malgrange's division theorem, e.g., \cite{AGV}]\label{Malgrange}
Let $f(t, \bx)$ be a real-valued smooth function defined near the origin of $\R^{n+1}$, where $t\in\R$ and $\bx\in\R^n$. Set $P(t) := \sum_{i=0}^d \lambda_i t^i\; (\lambda_i\in\R, \lambda_d\not=0)$ a polynomial of degree $d$. 
Then, there exist smooth functions $Q (t, \bx)$ and $r_i(\bx) \; (0\leq i\leq d-1)$ near the origin such that
$$f(t, \bx) = Q(t, \bx)P(t) + \sum_{i=0}^{d-1}r_i(\bx)t^i.$$
\end{thm}

\begin{exam}\upshape 
{(Typical singularity types of maps \cite{AGV, IS})} 
We say that two map-germs $f, g: (\R^m, 0) \to (\R^n, 0)$ are {\em $\A$-equivalent} (or right-left equivalent) if there exist diffeomorphism-germs 
$\sigma$ and $\tau$ of the source and the target spaces such that $g=\tau \circ f \circ \sigma^{-1}$. 
We are interested in the following typical types arising as singularities of Legendre maps. 
A map-germ $h: (\R^n, 0) \to (\R^{n+1}, 0)$ is of type {\em cuspidal edge} if $h$ is $\A$-equivalent to the germ at $0$ given by $(x_1, \cdots, x_n) \mapsto (x_1^3, x_1^2, x_2, \cdots, x_n)$. 
Also $h$ is of type {\em swallowtail} if it is $\A$-equivalent to 
$$(x_1, \cdots, x_n) \mapsto (3x_1^4+x_1^2x_2, 2x_1^2+x_1x_2, x_2, \cdots, x_n).$$
Those types are of corank one, i.e., $\rank dh_0=n-1$. 
We also denote by $h: U \to \R^{n+1}$ a representative of the germ with a neighborhood $U$ of $0$.  Obviously, in each case, the singular set $S(h)=\{p \in U, \rank dh_p \le n-1\}$ is a submanifold of codimension one, and the germ $h: (U, p) \to (\R^{n+1}, h(p))$ at every point  $p \in S(h)$ is of type cuspidal edge (except for $p=0$ in the latter case).  
 \end{exam}

Now we explain our setup. Since this is a local problem, 
we consider a local model $L \subset \R^{2n+1}$ with quasi-Hessian metric $h$ and coherent tangent bundles $(E, \Phi, \nabla^E)$ and $(E', \Phi', \nabla^{E'})$, and take a point $\bar{p}=(\bar{\bx}, \bar{\bp}, \bar{z}) \in L$. 
Below we will state our results only for the $m$-wavefront associated to $L$ around $\bar{p}$, but similar results for the $e$-wavefront hold just by exchanging $(\Phi', \nabla^{E'})$ and $(\Phi, \nabla^E)$. 

Of our interest is local geometry of singularities of the {\em $m$-Lagrange map}
$$\pi^m_1:=\Pi\circ\pi^m: L\to \R^n_{\bbp}$$ 
at singular points, where $\Pi:\R^n\times\R\to\R^n$ is the projection to the first factor. 
Suppose that $d\pi^m_1$ is degenerate at the point $\bar{p}$. 
In particular, we consider the case of corank one, i.e., 
$$\dim \ker d\pi^m_1(\bar{p})\,  (\, = \dim \ker \Phi' (\bar{p})\, )\, =1.$$
Then the singular point set $S(\pi^m_1) \subset L$ in a neighborhood of $\bar{p}$ is described by a {\em discriminant function} (in a system of local coordinates of $L$)  
$$\lambda (q):=\det\left[d\pi^m_1(q)\right] \quad (q \in L)$$ 
so that $S(\pi^m_1)=\{\lambda=0\}$.  
Also we take a non-zero vector field $X$ on $L$ around $\bar{p}$ which spans $\ker d\pi^m_{\bar{p}}$.

First, according to Izumiya-Saji's convenient criterion \cite{IS}, 
the $m$-Legendre map 
$$\pi^m=(\pi^m_1, z'): L \to \R^{n+1}=\R^n_{\bbp}\times \R_{z'}$$ 
has a singularity of type cuspidal edge at $\bar{p}$ if and only if it holds that
\begin{align*}
X\lambda(\bar{p}) \neq 0.
\end{align*}
Here, nothing about information geometry appears yet. 
Our first task is to interpret this condition in terms of a quasi-Hessian geometric quantity associated to $L$. 
Then we provide locally a normal form of the multi-valued (dual) potential for the $m$-wavefront using {\em affine flat coordinates}. 
This generalizes a result of an early work of Ekeland \cite[Proposition 2.7]{Ekeland}. 
We remark that the proof of \cite{Ekeland} seems to be incomplete, 
and our proof corrects it.

We fix a point $\bar{p}=(\bar{\bx}, \bar{\bp}, \bar{z}) \in L$ with 
$\bar{\bx}=(\bar{x}_1, \cdots, \bar{x}_n) \in \R_{\bbx}^n$ and  
$\bar{\bp}=(\bar{p}_1, \cdots, \bar{p}_n) \in \R_{\bbp}^n$ (bar means the fixed one, not a variable). 

\begin{thm}\label{m-wavefront1} 
Assume that $\rank \Phi'(\bar{p}) = n-1$ and 
\begin{align}\label{m-condition1}
\tau(\Phi(X)(\bar{p}), (\nabla^{E'}_{X}\Phi'(X))(\bar{p}))\neq 0
\end{align}
for any non-zero vector field $X$ on $L$ near $\bar{p}$ which spans $\ker \Phi'({\bar{p}})$.
Then, $\pi^m$ has a singularity of type cuspidal edge at $\bar{p}$, that is, 
we may retake affine flat coordinates 
$\bp=(p_1,\cdots,p_n)$ of $\R^n_{\bbp}$ and its dual $\bx=(x_1, \cdots, x_n)$ of $\R^n_{\bbx}$ 
such that there exist neighborhoods  
$\bar{p} \in V\subset L$,
$0 \in \U_1\subset\R$,
$(\bar{p}_2, \cdots, \bar{p}_{n})  \in \U_2\subset\R^{n-1}$ 
and functions 
$$\vp\in C^\infty(\U_1\times\U_2),\; \;\; k_1,k_2\in C^\infty(\U_2)$$ 
satisfying that  
the bi-valued dual potential for the $m$-wavefront $\pi^m(V)$ is given by 
\begin{align}\label{normal_1}
z'=z'_\pm(\bp)=k_2(p_2,\cdots,p_{n})+(p_1-k_1)\vp(\pm\sqrt{p_1-k_1}, p_2,\cdots,p_{n})
\end{align}
with $p_1\geq k_1(p_2,\cdots,p_{n})$,  
and 
the derivative of $\vp$ with respect to the first coordinate is non-zero.
\end{thm}

\proof 
From the assumption that $\rank \Phi'(\bar{p}) = n-1$, we may assume that $\ker \Phi'(\bar{p})$ is spanned by 
$(\frac{\rd}{\rd x_1})_{\bar{p}} + \bar{p}_1(\frac{\rd}{\rd z})_{\bar{p}}$
by taking an affine transformation of $\R^n_{\bbx}$, 
and thus by an affine Legendre equivalence of $\R^{2n+1}$. 
Then $(x_1, p_2,\cdots,p_{n})$ is a system of local coordinates of $L$ around $\bar{p}$. 
The singular set of the $m$-Lagrange map is $S(\pi^m_1) = \{\frac{\rd p_1}{\rd x_1}=0\}$ in this coordinates.

Let $g(x_1, p_2,\cdots,p_{n})=g(x_1, \bp_J)$ be a generating function of $L$ around $\bar{p}$ (i.e., $I=\{1\}, J=\{2, \cdots, n\}$). Then  
\begin{equation}\label{pz_gen}
p_1=\frac{\rd g}{\rd x_1}(x_1, \bp_J), \quad 
z'= x_1\frac{\rd g}{\rd x_1}(x_1, \bp_J) - g(x_1, \bp_J). 
\end{equation}
Note that $\frac{\rd z'}{\rd x_1}=x_1\frac{\rd p_1}{\rd x_1}=0$ on the singular set $S(\pi^m_1)$.
Take a vector field $X$ with $X=\frac{\rd}{\rd x_1}$ near $\bar{p}$ (which spans $\ker\Phi'(\bar{p})$), 
and take the flat sections $s_i$ of $E$ and $s^*_j$ of $E'$ as mentioned before. 
Note that as in the proof of \cite[Lem.3.22]{NT}, 
for any local section $\zeta'=\sum b_j s_j^*$ of $E'$ around $\bar{p}$, we see $\nabla_X^{E'}\zeta'=\sum X(b_j) s_j^*$. 
A direct computation shows that 
\begin{align*}
\Phi(X(\bar{p})) &=s_1 - \sum_{j=2}^{n}\frac{\rd^2 g}{\rd x_1\rd p_j}(\bar{p})s_j,\\
\Phi'(X(\bar{p})) &= \frac{\rd^2 g}{\rd x_1^2}(\bar{p})s^*_1 = \frac{\rd p_1}{\rd x_1}(\bar{p})s^*_1,\\
(\nabla^{E'}_{X}\Phi'(X))(\bar{p}) &= \frac{\rd^3 g}{\rd x_1^3}(\bar{p})s^*_1 = \frac{\rd^2 p_1}{\rd x_1^2}(\bar{p})s^*_1. 
\end{align*}
Hence, 
\begin{align}\label{tau_x}
\tau(\Phi(X(\bar{p})), (\nabla^{E'}_{X}\Phi'(X))(\bar{p})) &= \frac{1}{2}\frac{\rd^3 g}{\rd x_1^3}(\bar{p})
=\frac{1}{2}\frac{\rd^2 p_1}{\rd x_1^2}(\bar{p}).
\end{align}
From the second assumption, we have $\frac{\rd^2 p_1}{\rd x_1^2}(\bar{p})\neq 0$. 
Thus,  on the singular set $S(\pi^m_1)$ in an open neighborhood of $\bar{p}$, it holds that 
\begin{align}\label{dif}
\frac{\rd p_1}{\rd x_1}=0,\;\;
\frac{\rd^2 p_1}{\rd x_1^2} \neq 0, \;\;
\frac{\rd z'}{\rd x_1} = 0.  
\end{align}

For the first equation in $(\ref{dif})$, we use the implicit function theorem; 
there exist a neighborhood $\U_2$ of $(\bar{p}_2, \cdots, \bar{p}_n)\in \R^{n-1}$
and a function $f\in C^\infty(\U_2)$ whose graph $x_1=f(\bp_J)$ parametrizes   
the singular set $S(\pi^m_1)$ near $\bar{p}$. 

From Malgrange's division theorem, replacing $\U_2$ with smaller one if necessary, there exist functions $\varphi_1, \varphi_2, k_1, k_2, k'_1, k'_2$ such that 
\begin{align*}
p_1(x_1, \bp_J) = &(x_1-f)^2\varphi_1(x_1, \bp_J)+ (x_1-f)k'_1(\bp_J) + k_1(\bp_J),\\
z'(x_1, \bp_J) = &(x_1-f)^2\varphi_2(x_1, \bp_J)+ (x_1-f)k'_2(\bp_J) + k_2(\bp_J).
\end{align*}
Then, from (\ref{dif}), we have $k'_1=k'_2=0$ and $\varphi_1\neq 0$ on the singular set.

Assume that $\varphi_1 > 0$ (if $\vp_1 < 0$, change $p_1$ to $-p_1$). 
Define a local coordinate change $\Psi:(x_1, \bp_J)\mapsto(y_1, \bp_J)$ with $y_1:=(x_1-f)\sqrt{\varphi_1}$ and put $\vp:=(\vp_2/\vp_1)\circ\Psi^{-1}$. 
Then, $p_1 - k_1=y_1^2 \geq 0$ and 
\begin{align*}
z'-k_2 &= y_1^2\vp(y_1, \bp_J)= (p_1-k_1)\vp(\pm\sqrt{p_1-k_1}, \bp_J).
\end{align*} 
Finally we take $\U_1$ so that the image of $\Psi$ is $\U_1\times\U_2$ and put $V\subset L$ to be the corresponding neighborhood of $\bar{p}$.
Thus we get the normal forms (\ref{normal_1}). 

From the equality above, we have $\frac{\rd^3 z'}{\rd y_1^3}(\bar{p})=6\frac{\rd\vp}{\rd y_1}(\bar{p})$, 
and from (\ref{pz_gen}),
\begin{align*}
z'(y_1, \bp_J) &= x_1p_1 - \tilde{g}(y_1, \bp_J)= x_1(y_1^2 + k_1) - \tilde{g}(y_1, \bp_J),
\end{align*}
where $\tilde{g}:=g\circ\Psi^{-1}$. 
Hence, $\frac{\rd^3 z'}{\rd y_1^3}=2y_1\frac{\rd^2 x_1}{\rd y_1^2} + 4\frac{\rd x_1}{\rd y_1}$. 
Since $y_1(\bar{p})=0$ and $\Psi$ is a coordinate change, it holds that
\begin{align*}
\frac{\rd\vp}{\rd y_1}(\bar{p})=\frac{1}{6}\frac{\rd^3 z'}{\rd y_1^3}(\bar{p}) = \frac{2}{3}\frac{\rd x_1}{\rd y_1}(\bar{p})\neq 0.
\end{align*}
This completes the proof. 
\qed

\begin{exam}\label{A2}\upshape
{\bf (Cuspidal edge).} 
Consider the case of $n=2$ and a generating function given by  
$$g(x_1, p_2) = \frac{1}{3}x_1^3 - \frac{1}{2}p_2^2.$$
Then the $e$-wavefront $W_e(L)$ is smooth and the $m$-wavefront $W_m(L)$ is the graph of the bi-valued dual potential 
$$z'=x_1p_1-g=\pm\frac{2}{3}p_1^{3/2}+\frac{1}{2}p_2^2$$  
defined on $p_1\ge 0$ and branched along $p_1=0$.  
\end{exam}

\begin{rem}\upshape
From the equality (\ref{tau_x}) in the proof above, we have
\begin{align}\label{tau_xx}
\tau(\Phi(X(\bar{p})), (\nabla^{E'}_{X}\nabla^{E'}_{X}\Phi'(X))(\bar{p})) = \frac{1}{2}\frac{\rd^4 g}{\rd x_1^4}(\bar{p})
= \frac{1}{2}(XX\lambda)(\bar{p}),
\end{align}
where $X=\frac{\rd}{\rd x_1}$ and $\lambda=\frac{\rd p_1}{\rd x_1}$. 
Hence, using Izumiya-Saji's criterion \cite{IS}, 
the germ of the $m$-wavefront $\pi^m$ at $\bar{p}$ is diffeomorphic to the swallowtail singularity if and only if the following (in)equalities at $\bar{p}$ hold: $d\lambda(\bar{p})\neq 0$ and  
\begin{align*}
  \tau(\Phi(X)(\bar{p}), (\nabla^{E'}_{X}\Phi'(X))(\bar{p})) = 0,\\
\tau(\Phi(X)(\bar{p}), (\nabla^{E'}_{X}\nabla^{E'}_{X}\Phi'(X))(\bar{p})) \neq 0. 
\end{align*}
Although we may write down a normal form of the multi-valued potential in affine coordinates, we do not execute this due to its complexity. 
As an example, consider a generating function 
$$g(x_1, p_2)=\frac{1}{4} x_1^4 + \frac{1}{2} p_2 x_1^2.$$
The $m$-wavefront forms the swallowtail surface (cf. \cite{AGV}). 
\end{rem}

On the other hand, in practical applications of information geometry, it is quite natural to assume that a quasi-Hessian metric is positive semi-definite (e.g., the Fisher-Rao metric is so). 
Among such cases, let us consider the following simplest setting.
Suppose that $h$ is positive semi-definite and $\rank \Phi'(\bar{p})=n-1$.
In this case, the derivative of a discriminant function $\lambda$ for $\pi^m_1$ is zero on the singular set.
Indeed, from Proposition \ref{quasi-Hessian_local}, for a generating function $g(x_1,\bp_J)$, the quasi-Hessian metric is expressed by
\begin{align*}
h =
\begin{bmatrix}
\frac{\rd^2 g}{\rd x_1^2} & 0\\
0 & -g_{JJ}
\end{bmatrix},
\end{align*}
where $J=\{2,\cdots,n\}$ and $g_{JJ} = [\frac{\rd^2 g}{\rd p_i \rd p_j}]_{i,j\in J}$. Since $h$ is positive semi-definite, we have $\frac{\rd^2 g}{\rd x_1^2}=\lambda \geq 0$. Hence $d\lambda = 0$ on the singular set $\{\lambda=0\}$.
This means that the implicit function theorem for $\lambda=0$ does not hold.
Even though, it is reasonable, as a simplest case, to consider that the singular set is a submanifold of codimension one in $L$.
The following theorem describes the affine normal form appearing in this case, and also characterizes it as in the context of the {\em classification of minimal points of a function} due to V. A. Vasil'ev \cite{Vasil'ev}.

\begin{thm}\label{m-wavefront2}
Assume that $\rank \Phi'({\bar{p}}) = n-1$ and 
\begin{align}
\tau(\Phi(X)(q), (\nabla^{E'}_{X}\Phi'(X))(q)) = 0, \label{m-condition2_1}\\
\tau(\Phi(X)(\bar{p}), (\nabla^{E'}_{X}\nabla^{E'}_{X}\Phi'(X))(\bar{p})) \neq 0 \label{m-condition2_2}
\end{align}
for any non-zero vector field $X$ on $L$ which spans $\ker \Phi'(q)$ at any $q$ in the singular set $S(\pi^m_1)$ near $\bar{p}$. 
Suppose that $S(\pi^m_1)$ is a submanifold of codimension one in $L$.
Then, we may retake affine flat coordinates $\bp$ and its dual coordinates $\bx$ such that 
there exist neighborhoods  $\bar{p}\in V\subset L$, 
$0\in\U_1\subset\R$,
$(\bar{p}_2, \cdots, \bar{p}_n) \in\U_2\subset\R^{n-1}$ and functions 
$$\vp\in C^\infty(\U_1\times\U_2), \;\;\; k_1,k_2\in C^\infty(\U_2)$$ 
satisfying that 
the dual potential for the $m$-wavefront $\pi^m(V)$ is given by  
\begin{align}\label{normal_2}
z'=z'(\bp)=k_2(p_2,\cdots,p_n)+(p_1-k_1)\vp((p_1-k_1)^{1/3}, p_2,\cdots,p_n),
\end{align}
where the derivative of $\vp$ with respect to the first coordinate is non-zero.
\end{thm}

\begin{rem}\upshape
As seen before, in the case where a quasi-Hessian metric $h$ is positive semi-definite, the derivative of a discriminant function for $\pi^m_1$ is zero on the singular set, and then the condition (\ref{m-condition2_1}) above is automatically satisfied.
\end{rem}

\proof
We use the notation and the settings in the proof of Theorem \ref{m-wavefront1}. 
Note that $\frac{\rd}{\rd x_1}$ spans $\ker \Phi'(q)$ at any $q$ in the singular set $S(\pi^m_1)=\{\frac{\rd p_1}{\rd x_1}=0\}$ near $\bar{p}$.
In a similar way to the proof of Theorem \ref{m-wavefront1}, from the condition (\ref{m-condition2_1}) and (\ref{m-condition2_2}), we see that 
\begin{align}\label{dif2}
\frac{\rd p_1}{\rd x_1}=0,\; \frac{\rd^2 p_1}{\rd x_1^2}=0,\; \frac{\rd^3 p_1}{\rd x_1^3}\neq0,\; \frac{\rd z'}{\rd x_1}=0,\; \frac{\rd^2 z'}{\rd x_1^2}=0
\end{align}
hold on the singular set $S(\pi_1^m)$ in an open neighborhood of $\bar{p}$. 
For the second equation in ($\ref{dif2}$), we use the implicit function theorem; there exist a neighborhood $\U_2\subset \R^{n-1}$ of $(\bar{p}_2, \cdots, \bar{p}_n)$, and a function $f\in C^{\infty}(\U_2)$ satisfying $\frac{\rd^2 p_1}{\rd x_1^2}(f(\bp_J),\bp_J)=0$ on $\U_2$. 
This graph coincides with $S(\pi_1^m)$ near $\bar{p}$, since $S(\pi_1^m)$ is assumed to be a codimension one submanifold of $L$. 

From Malgrange's division theorem, there exist functions $\vp_1$, $\vp_2$, $k_1$, $k_2$, $k'_1$, $k'_2$, $k''_1$, $k''_2$ such that 
\begin{align*}
p_1(x_1,\bp_J) = &(x_1-f)^3\varphi_1(x_1,\bp_J)+ (x_1-f)^2k'_1(\bp_J) + (x_1-f)k''_1(\bp_J) + k_1(\bp_J),\\
z'(x_1,\bp_J) = &(x_1-f)^3\varphi_2(x_1,\bp_J)+ (x_1-f)^2k'_2(\bp_J) + (x_1-f)k''_2(\bp_J) + k_2(\bp_J).
\end{align*}
Then from (\ref{dif2}), we have $k'_1 =k''_1=k'_2 =k''_2=0$ and $\vp_1\neq 0$ on the singular set. 
Namely, 
\begin{align*}
p_1(x_1,\bp_J) &= (x_1-f)^3\varphi_1(x_1,\bp_J)+ k_1(\bp_J),\\
z'(x_1,\bp_J) &= (x_1-f)^3\varphi_2(x_1,\bp_J)+k_2(\bp_J).
\end{align*}
Set $\Psi:(x_1,\bp_J)\mapsto(y_1,\bp_J)$ with $y_1:=(x_1-f)(\varphi_1)^{1/3}$ and put $\vp:=(\vp_2/\vp_1)\circ\Psi^{-1}$.
Then
$$z'-k_2 = y_1^3\vp(y_1,\bp_J) = (p_1-k_1)\vp((p_1-k_1)^{1/3}, \bp_J).$$
Here we take $\U_1$ and $V\subset L$ suitably. Thus we get the normal forms (\ref{normal_2}). 

Also we see that $\frac{\rd \vp}{\rd y_1}(\bar{p})\neq 0$. Indeed, from the equality above, we have
$\frac{\rd^4 z'}{\rd y_1^4}(\bar{p})=24\frac{\rd\vp}{\rd y_1}(\bar{p})$.
Letting $g(x_1, \bp_J)$ be a generating function, we have
\begin{align*}
z'(y_1,\bp_J) &= x_1p_1 - \tilde{g}(y_1,\bp_J)= x_1(y_1^3 + k_1) - \tilde{g}(y_1,\bp_J),
\end{align*}
where $\tilde{g}:=g\circ\Psi^{-1}$, and then $\frac{\rd^4 z'}{\rd y_1^4}=3y_1^2\frac{\rd^3 x_1}{\rd y_1^3}+18y_1\frac{\rd^2 x_1}{\rd y_1^2}+18\frac{\rd x_1}{\rd y_1}$. 
Since $y_1(\bar{p})=0$ and $\Psi$ is a coordinate change, it holds that
\begin{align*}
\frac{\rd\vp}{\rd y_1}(\bar{p}) = \frac{1}{24}\frac{\rd^4 z'}{\rd y_1^4}(\bar{p}) = \frac{3}{4}\frac{\rd x_1}{\rd y_1}(\bar{p})\neq 0.
\end{align*}
This completes the proof. 
\qed

\begin{exam}\label{semi-definite}\upshape
{\bf  ($A_3$-singularity \cite{Vasil'ev}).} Consider a generating function 
$$g(x_1, p_2) = \frac{1}{4}x_1^4 - \frac{1}{2}p_2^2.$$ 
The $e$-wavefront $W_e(L)$ is smooth and $m$-wavefront $W_m(L)$ is then the graph of the dual potential: 
$$z'=x_1p_1-g=\frac{3}{4}p_1^{4/3}+\frac{1}{2}p_2^2.$$ 
The $e$-caustics is empty, while the $m$-caustics appear along $p_1=0$. 
\end{exam}

\section{ Expressions of the criteria and contrast functions}\label{sec:expression of criteria}
As seen in \S \ref{sec:2}, a quasi-Hessian manifold $M$ is endowed with a quasi-Hessian metric $h$ and a symmetric cubic tensor $C$, and these tensors are locally written by the second and third derivative of a generating function. 
In this section, we discuss the relationship between such tensors and our Theorems \ref{m-wavefront1} and \ref{m-wavefront2}. 

The quantity $\tau(\Phi(X), \nabla^{E'}_X\Phi'(X))$ in the condition (\ref{m-condition1}) of Theorem \ref{m-wavefront1} can be rewritten by using $h$ and $C$.

\begin{prop}\label{criterion tensor}
For any vector field $X$ on $M$, it holds that
$$
\tau(\Phi(X), \nabla^{E'}_X\Phi'(X)) = \frac{1}{4}\left(C(X, X, X) + Xh(X, X)\right).
$$
\end{prop}
\proof
From the definition \ref{quasi-Hessian metric} and \ref{cubic tensor}, we have
\begin{align*}
h(X, X) &= 2\tau(\Phi(X), \Phi'(X)),\\
C(X, X, X) &= 2\tau(\Phi(X), \nabla^{E'}_X\Phi'(X)) - 2\tau(\nabla^E_X\Phi(X), \Phi'(X)).
\end{align*}
Since 
$Xh(X, X) = 2\tau(\nabla^E_X\Phi(X), \Phi'(X)) + 2\tau(\Phi(X), \nabla^{E'}_X\Phi'(X))$
(see Lemma $3.22$ in \cite{NT}), the equality follows. \qed

\

One may expect to express the quantity $\tau(\Phi(X), \nabla^{E'}_X\nabla^{E'}_X\Phi'(X))$ appearing in the condition (\ref{m-condition2_2}) of Theorem \ref{m-wavefront2} by using some tensors $h$, $C$ and $X$, like as Proposition \ref{criterion tensor}.
However, it turns out to be impossible -- a natural candidate is the sum of two quantities, 
\begin{align*}
XC(X,X,X) = &2\{\tau(\Phi(X), \nabla^{E'}_X\nabla^{E'}_X\Phi'(X)) - \tau(\nabla^{E}_X\nabla^{E}_X\Phi(X), \Phi'(X))\},\\
XXh(X, X) = &2\{\tau(\Phi(X), \nabla^{E'}_X\nabla^{E'}_X\Phi'(X)) + \tau(\nabla^{E}_X\nabla^{E}_X\Phi(X), \Phi'(X))\}\\
&+ 4\tau(\nabla^E_X\Phi(X), \nabla^{E'}_X\Phi'(X))\nonumber,
\end{align*}
but $\tau(\nabla^E_X\Phi(X), \nabla^{E'}_X\Phi'(X))$ does not vanish in general. 
Nevertheless, there might exist another fourth-order tensor on $M$ which expresses the fourth derivatives of a generating function locally. 
In fact, this question touches an important aspect of geometry of statistical manifolds \cite{Eguchi, CMP, Matsuzoe1999}, that is explained below.

A {\em statistical manifold} is a manifold $N$ endowed with a (pseudo-)Riemannian metric $h$ and a symmetric cubic tensor $C$, 
and a dually flat manifold is a most particular such one. 
According to Eguchi \cite{Eguchi},  a statistical manifold $N$ is characterized by some `asymmetric distance function' $\rho: N \times N \to \R$. 
Given vector fields $X_1, \cdots, X_k, Y_1,\cdots, Y_l$ on $N$, 
we set a function 
$$\rho[X_1\cdots X_k|Y_1\cdots Y_l]:N\to\R$$
to be defined by
\begin{align}\label{contrast}
\rho[X_1\cdots X_k|Y_1\cdots Y_l](r) = (X_1)_p\cdots(X_k)_p(Y_1)_q\cdots(Y_l)_q(\rho(p,q))|_{p=q=r}. 
\end{align}
We say that $\rho:N\times N\to\R$ is a {\em contrast function} if it holds that
\begin{itemize}
\item[(i)] $\rho[-|-](r)=\rho(r,r) = 0$,
\item[(ii)] $\rho[X|-](r) = \rho[-|X](r) = 0$,
\item[(iii)] $-\rho[X|Y]$ is a pseudo-Riemannian metric on $N$
\end{itemize} 
for $r\in N$ and vector fields $X, Y$ on $N$. 
Eguchi has investigated the geometry derived from contrast functions \cite{Eguchi}; e.g., 
a contrast function directly yields geometric quantities such as metric $h$, cubic tensor $C$ and so on. 
Conversely, it is shown by Matumoto \cite{Matumoto} that given $h$ and $C$, one can construct a contrast function which recovers them.

For a dually flat manifold, take a local convex potential, then the so-called {\em Bregman divergence} gives a contrast function.
In our case, local potential may not be convex (nor single-valued), but in the following way we have the {\em canonical divergence}, which gives a {\em weak contrast function} \cite{NT}, i.e., only the condition (i) and (ii) above are satisfied.

For a Legendre submanifold $L\subset\R^{2n+1}$, the canonical divergence $D:L\times L\to\R$ is defined by
\begin{align*}
  D(p, q) = z(p) + z'(q) - \bx(p)^T\bp(q).
\end{align*}
Since the canonical divergence $D$ is invariant under affine Legendre equivalences, it is defined on a quasi-Hessian manifold $M$ \cite{NT}. Let denote $D_M$ be the canonical divergence on $M$.

The canonical divergence $D_M$ restores a quasi-Hessian metric $h$ and a canonical cubic tensor $C$ (Theorem $4.10$ in \cite{NT}):
\begin{align*}
h(X,Y)=-D_M[X|Y], \;\; C(X, Y, Z) = -D_M[Z|XY] + D_M[XY|Z].
\end{align*}
In this way, by using functions of the form (\ref{contrast}), tensors $h$ and $C$ are made up of a contrast function.

Recently, it is announced in \cite{CMP} that for a statistical manifold,
there are only two fourth-order tensors made up of a contrast function, and especially, they vanish for the Bregman divergence associated to a dually flat manifold.
According to this, a crucial consequence is that 
it is impossible to construct a fourth-order tensor which locally expresses the fourth derivative of a generating function $g$ on a dually flat manifold (quasi-Hessian manifold).

Instead of using tensors, we see that the fourth derivative of $g$ can be directly written by using the canonical divergence $D_M$ on a quasi-Hessian manifold $M$. 
That follows from a more general formula (\ref{dm2}) in the following theorem obtained by a direct computation. 
\begin{thm}\label{equ_div}
Let $X, Y, Z, W$ be vector fields on $M$, then it holds that 
\begin{align}
\tau(\Phi(Z), \nabla^{E'}_Y\Phi'(W)) &= -\frac{1}{2}D_M[Z|YW],\label{dm1}\\
\tau(\Phi(Z), \nabla^{E'}_X\nabla^{E'}_Y\Phi'(W)) &= -\frac{1}{2}D_M[Z|XYW]\label{dm2}.
\end{align}
\end{thm}
\proof 
Take a local model $L\subset\R^{2n+1}$ of $M$, and let $g(\bx_I, \bp_J)$ be a generating function of $L$, where $(\bx_I, \bp_J)$ are local coordinates around $p\in L$. Then 
$$\bx_J(q)=-\frac{\rd g}{\rd \bp_J}(q), \;\; \bp_I(q)=\frac{\rd g}{\rd \bx_I}(q), \;\; z(q)=\bx_J(q)^T\bp_J(q)+g(\bx_I(q), \bp_J(q))$$ 
and 
\begin{align}\label{D}
D_M(p,q) = g(p) - g(q) + \bx_J(p)^T(\bp_J(p)-\bp_J(q)) + \bp_I(q)^T(\bx_I(q)-\bx_I(p))
\end{align}
for $q\in L$ close to $p$, where $g(q):=g(\bx_I(q), \bp_J(q))$.
Let $\rd_i$ denote $\frac{\rd}{\rd x_i}$ if $i\in I$ and $\frac{\rd}{\rd p_i}$ if $i\in J$.

Since $h(Z, W)=-D_M[Z|W]$, it suffices to show that the equality (\ref{dm1}) holds for $Y=\rd_t, Z=\rd_u, W=\rd_v$.
Letting $s^*_i$ be the flat sections of $E'$ $(1\leq i\leq n)$, we have
\begin{align*}
\nabla^{E'}_{\rd_t}\Phi'(\rd_v) = \sum_{i\in I}\rd_t\rd_v\rd_i g s^*_i.
\end{align*}
Hence 
\begin{align*}
2\tau(\Phi(\rd_u), \nabla^{E'}_{\rd_t}\Phi'(\rd_v)) =
\left\{
\begin{array}{ll}
\rd_t\rd_u\rd_v g &(u\in I),\\
0 &(u\in J).
\end{array}
\right.
\end{align*}
On the other hand, by differentiating the right-hand side of (\ref{D}) (see \cite{NT}), we have
\begin{align*}
(\rd_u)_p(\rd_t)_q(\rd_v)_qD_M(p, q) = -\rd_t\rd_v\bp_I(q)^T(\rd_u\bx_I(p)).
\end{align*}
Hence
\begin{align*}
D_M[(\rd_u)|(\rd_t)(\rd_v)] = 
\left\{
\begin{array}{ll}
-\rd_t\rd_u\rd_v g &(u\in I),\\
0 &(u\in J).
\end{array}
\right.
\end{align*}
This means that the equality (\ref{dm1}) holds. Put $X=\rd_s$. Similarly, 
\begin{align*}
2\tau(\Phi(\rd_u), \nabla^{E'}_{\rd_s}\nabla^{E'}_{\rd_t}\Phi'(\rd_v)) =
\left\{
\begin{array}{ll}
\rd_s\rd_t\rd_u\rd_v g &(u\in I),\\
0 &(u\in J).
\end{array}
\right.
\end{align*}
and 
\begin{align*}
D_M[(\rd_u)|(\rd_s)(\rd_t)(\rd_v)]=
\left\{
\begin{array}{ll}
-\rd_s\rd_t\rd_u\rd_v g &(u\in I),\\
0 &(u\in J).
\end{array}
\right.
\end{align*}
Namely, the equality (\ref{dm2}) holds.
\qed

\

Thus our geometric criteria for singularity types are expressed by using the canonical divergence: 
\begin{cor}\label{DXXXX}
For any vector field $X$ on $M$, it holds that
\begin{align}
\tau(\Phi(X), \nabla^{E'}_X\Phi'(X)) &= -\frac{1}{2}D_M[X|XX],\label{div1}\\
\tau(\Phi(X), \nabla^{E'}_X\nabla^{E'}_X\Phi'(X)) &= -\frac{1}{2}D_M[X|XXX]\label{div2}.
\end{align}
\end{cor}

\begin{rem}\upshape 
In \cite{Eguchi}, Eguchi has explored relations among forth-order tensors and derivatives of contrast functions, and especially introduced 
a special fourth-order tensor $B^*$ (also named the {\em Bartlett tensor}\, in \cite{Matsuzoe1999}). 
It can be shown that the equality (\ref{dm2}) in Theorem \ref{equ_div} (in case that $M$ is a dually flat manifold) is equivalent to $B^*=0$ ($B^*$-free). In this sense,  the formula (\ref{div2}) in Corollary \ref{DXXXX} reflects such fourth-order geometry of statistical manifolds due to Eguchi. 
\end{rem}


\end{document}